\documentclass{article}

\newtheorem{fed}{\textbf{Definition}}[section]
\newtheorem{thm}[fed]{\textbf{Theorem}}
\newtheorem{lemma}[fed]{\textbf{Lemma}}

\newtheorem{prop}[fed]{\textbf{Proposition}}

\usepackage{amscd,amssymb,bbm,graphicx,epsfig,psfrag,epic,eepic,latexsym}
\usepackage{amsmath}
\usepackage[arrow,matrix,curve]{xy}
\usepackage{mathrsfs}
\usepackage[dvips]{color}

\begin{document}
\title{Morse homology on noncompact manifolds}
\author{Kai Cieliebak, Urs Frauenfelder}
\maketitle

\begin{abstract}
Given a Morse function on a manifold whose moduli spaces
of gradient flow lines for each action window
are compact up to breaking one gets a bidirect system of chain
complexes. There are different possibilities to take
limits of such a bidirect system. We discuss in this
note the relation between these different limits.
\end{abstract}

\tableofcontents

\section[Introduction]{Introduction}\label{intro}

In this note we assume that we have a Morse function $f$ on
a finite dimensional (possibly noncompact) Riemannian manifold $(M,g)$
with the property that the moduli spaces of gradient flow lines
in fixed action windows are compact up to breaking. Hence for an action
window $[a,b] \subset \mathbb{R}$ we can define Morse homology groups
$$HM^{[a,b]}_*=HM^{[a,b]}_*(f,g).$$
For our purposes the following notation turns out to be useful
$$HM_a^b:=HM^{[a,b]}_*.$$
So the reader should be aware that the subscript for our homology
groups does not refer to the grading but to the lower end of the action
window. We actually suppress the reference to the grading since it
plays a minor role in our discussion.

There are now different limits one can take from these homology groups.
One possibility was carried out by H.\,Hofer and D.\,Salamon in
\cite{hofer-salamon}. They take a Novikov completion of the chain
complex on which they get a well-defined boundary homomorphism.
We denote by $HM$ the homology of this complex. Other
possibilities are to take direct and inverse limits of
the homology groups $HM^b_a$. Hence we abbreviate
$$\overline{HM}=\lim_{\substack{\longrightarrow\\ b \to
\infty}}\lim_{\substack{\longleftarrow\\
a \to -\infty}}HM_a^b$$
and
$$\underline{HM}=\lim_{\substack{\longleftarrow\\ a \to
-\infty}}\lim_{\substack{\longrightarrow\\
b \to \infty}}HM_a^b.$$
The aim of this note is to study the relation between
these three homology groups. We remark that there are canonical maps
$\kappa \colon \overline{HM} \to \underline{HM}$,
$\overline{\rho} \colon HM \to \overline{HM}$, and
$\underline{\rho} \colon HM \to \underline{HM}$ whose definition
we recall later. We summarize them into the diagram
\begin{equation}\label{dia}
\begin{xy}
 \xymatrix{HM \ar[r]^{\overline{\rho}}
\ar[rd]_/-0.8em/{\underline{\rho}}& \overline{HM} \ar[d]^{\kappa}\\
 &\underline{HM}
}
\end{xy}
\end{equation}
Our main result is the following.
\\ \\
\textbf{Theorem A: }\emph{Assume that the Morse homology
groups are taken with field coefficients. Then the diagram
(\ref{dia}) is commutative, $\overline{\rho}$ is an isomorphism, and
$\kappa$ and therefore also $\underline{\rho}$ are surjective.}
\\ \\
\emph{Remark 1: } Theorem A might fail if one uses integer coefficients
instead of field coefficients. We provide an example in the appendix.
\\ \\
\emph{Remark 2: } Although we state Theorem A only for finite
dimensional manifolds, it can be carried over to the semi-infinite
dimensional case of Floer homology. The difficulty with giving
a precise statement lies in the fact that up to now there
is no precise definition what a Floer homology in general
actually is. We hope to modify this unsatisfactory situation
in the near future by using the newly established theory
of H.\,Hofer, K.\,Wysocky, and E.\,Zehnder about scale structures
\cite{hofer-wysocky-zehnder}
to interpret Floer homology as Morse homology on scale manifolds.
Alternatively, one can also give an axiomatized description of
Morse homology for which Theorem A continues to hold. We explain
that in Section~\ref{axiom}.
\\ \\
\emph{Remark 3: } We became interested in the relation of the
different Morse homologies via Rabinowitz Floer homology.
We defined in \cite{cieliebak-frauenfelder} Rabinowitz Floer homology
as the Morse homology of the Rabinowitz action functional by taking
Novikov sums as in \cite{hofer-salamon}. On the other hand,
in a joint work with A.\,Oancea \cite{cieliebak-frauenfelder-oancea}
we are proving that Rabinowitz Floer homology is isomorphic
to a variant of symplectic homology. To establish this isomorphism
we need to work with $\overline{HM}$. Therefore it became
important for us to know if $\overline{\rho}$ is an isomorphism or
not.
\\ \\
\emph{Remark 4: } In Floer homology the homology groups
$\overline{HM}$ were successfully applied by K.\,Ono in his proof of
the Arnold conjecture for weakly monotone symplectic
manifolds~\cite{ono}. In this paper K.\,Ono raises the question if
$\overline{\rho}$ is an isomorphism. In the case of Floer homology
for weakly monotone symplectic manifolds
it was later shown by S.\,Piunikhin, D.\,Salamon, and M.\,Schwarz
that the homology groups $HM$ and $\overline{HM}$ coincide by
direct computation. Theorem A gives an algebraic explanation
for this fact.
\\ \\
The following example shows that $\kappa$ and
therefore $\underline{\rho}$ do not need
to be injective.
\\ \\
\emph{Example: }Let $M=\bigsqcup_{n=1}^\infty R_n$
where each $R_n \cong \mathbb{R}$ and for each
$n \in \mathbb{N}$ the Morse function $f|_{R_n}$ has
one single maximum $\overline{c}_n$ and one single minimum
$\underline{c}_n$ with
$$f(\overline{c}_n)=n,\quad f(\underline{c}_n)=-n.$$
It follows that there is precisely one gradient flow line
from $\overline{c}_n$ to $\underline{c}_n$. Taking
Morse homology with coefficients in the abelian
group $\Gamma$ we obtain for
$(a,b) \in \mathbb{R}^2$
$$HM^b_a=\bigg(\bigoplus_{b<n \leq -a} \Gamma
\cdot \underline{c}_n\bigg)\oplus
\bigg(\bigoplus_{-a<n \leq b} \Gamma \cdot\overline{c}_n\bigg)
$$
We conclude that
$$HM^b=\underleftarrow{\lim}HM^b_a=\prod_{n > b}
\Gamma \cdot \underline{c}_n$$
and
$$HM_a=\underrightarrow{\lim} HM^b_a=\bigoplus_{n>-a}\Gamma \cdot
\overline{c}_n.$$
We get
$$\underline{HM}=\underleftarrow{\lim}HM_a=0,
\quad \overline{HM}=\underrightarrow{\lim}HM^b \neq 0$$
which shows that $\kappa$ does not need to be injective.

\section[Morse homology]{Morse homology}

\subsection[Morse tuples]{Morse tuples}

In this section we introduce the notion of a Morse tuple on
a (not necessarily compact) finite dimensional manifold $M$. A Morse
tuple $(f,g)$ consists of a Morse function $f$ on $M$
and a Riemannian metric $g$ meeting a transversality and
a compactness condition which ensure that Morse homology
for each action window can be defined as usual,
see \cite{schwarz}. We then proceed by
explaining the maps which connect the Morse homology groups
for different action windows.
\smallskip

If $(M,g)$ is a Riemannian manifold and $f \in C^\infty(M)$
is a Morse function on $M$ we denote by $\nabla f$ the gradient of
$f$ with respect to
the metric $g$. A gradient flow line $x \in C^\infty(\mathbb{R},M)$
is a solution of the ordinary differential equation
\begin{equation}\label{grad}
\partial_s x(s)=\nabla f(x(s)), \quad s \in \mathbb{R}.
\end{equation}
We denote by $|| \cdot||$ the norm on $TM$ induced from the metric
$g$. The energy of any $x \in C^\infty(\mathbb{R},M)$
not necessarily satisfying (\ref{grad}) is given by
$$E(x)=\int_{-\infty}^\infty||\partial_s x||^2 ds.$$
If $x$ is a gradient flow line then its energy equals
$$E(x)=\limsup_{s \in \mathbb{R}} f(x(s))-\liminf_{s \in \mathbb{R}}
f(x(s))=\lim_{s \to \infty}f(x(s))-\lim_{s \to -\infty} f(x(s)).$$
In particular, if $x(s)$ converges to critical points $x^\pm$
of $f$ as $s$ goes to $\pm \infty$ we obtain
$$E(x)=f(x^+)-f(x^-).$$
We abbreviate
$$\mathcal{G}=\big\{x \in C^\infty(\mathbb{R},M):
\,\,x\,\,\textrm{solves (\ref{grad})},\,\,E(x)<\infty\big\}$$
the moduli space of all finite energy flow lines of $\nabla f$.
For a two dimensional vector $(a,b) \in \mathbb{R}^2$ we denote
$$\mathcal{G}_a^b=
\big\{x \in \mathcal{G}: a \leq f(x(s)) \leq b\text{ for all }
s \in \mathbb{R}\big\}.$$
Note that $\mathbb{R}$ acts on $\mathcal{G}$ by time shift
$$r_* x(s)=x(s+r), \quad x\in \mathcal{G},\,\,s,r \in \mathbb{R}.$$
This action is semifree in the sense that in the complement of
its fixed points it acts freely. We abbreviate
$$\mathcal{C}=\mathrm{Fix}(\mathbb{R}) \subset \mathcal{G}.$$
The fixed point set $\mathcal{C}$ can naturally be identified
with the set of critical points $\mathrm{crit}(f)$
via the evaluation map
$$\mathrm{ev} \colon \mathcal{G} \to M, \quad x \mapsto x(0).$$
Moreover, we endow the set
$\mathcal{C}$ with the structure of a graded
set where the grading is given by the Morse index.
For $(a,b) \in \mathbb{R}^2$ we further denote
$$\mathcal{C}_a^b=\mathcal{C} \cap \mathcal{G}_a^b.$$
Note that $\mathcal{C}_a^b$ corresponds to the critical points
of $f$ in the action window $[a,b]$. Again this set
is graded by the Morse index. We further remark that
$\mathcal{C}_a^b$ only depends on the Morse function $f$ and
not on the metric $g$.
Our first hypothesis is the following compactness assumption.
To state it we endow the space $C^\infty(\mathbb{R},M)$ with the
$C^\infty_{\mathrm{loc}}$-topology.
\begin{description}
 \item[(H1)] For all $(a,b) \in \mathbb{R}$ the set
  $\mathcal{G}_a^b$ is a compact subset of
  $C^\infty(\mathbb{R},M)$.
\end{description}
Before stating our second hypothesis we show in the following
lemma that hypothesis (H1) implies that each finite
energy gradient flow line converges asymptotically to critical
points.
\begin{lemma}\label{limit}
Assume hypothesis (H1). If $x \in \mathcal{G}$, then there exists
$x^\pm \in \mathcal{C}$ such that
$$\lim_{r \to \pm \infty}r_*x=x^\pm.$$
\end{lemma}
\textbf{Proof: }Choose $a,b \in \mathbb{R}$ such that
$x \in \mathcal{G}^b_a$. For $\nu \in \mathbb{N}$ consider the sequence
$$x_\nu=\nu_* x$$
of gradient flow lines. Since $\mathcal{G}^b_a$ is
$\mathbb{R}$-invariant it follows that
$$x_\nu \in \mathcal{G}^b_a, \quad \nu \in \mathbb{N}.$$
By hypothesis (H1) it follows
that there exists a subsequence $\nu_j$ and $x^+ \in \mathcal{G}^b_a$
such that $x_{\nu_j}$ converges to $x^+$ in the
$C^\infty_{\mathrm{loc}}$-topology
as $j$ goes to infinity. It remains to show that $x^+$ is a constant
gradient flow line, hence a critical point. Fix $s>0$.
To see that $x^+$ is constant we have to show that $f(x^+(0))=f(x^+(s))$.
We argue by contradiction and assume that there exists $\epsilon>0$ satisfying
$$f(x^+(s))-f(x^+(0))=\epsilon.$$
Since $x_{\nu_j}$ converges to $x^+$ in the $C^\infty_{\mathrm{loc}}$-topology
there exists $j_0$ such that for every $j \geq j_0$ the inequality
$$f(x_{\nu_j}(s))-f(x_{\nu_j}(0)) \geq \frac{\epsilon}{2}$$
holds. By definition of $x_{\nu_j}$ this means
$$f(x(\nu_j))-f(x(\nu_j+s)) \geq \frac{\epsilon}{2}.$$
For $\ell \in \mathbb{N}$ we define recursively
$$j_\ell=\min\big\{j: \nu_{j_{\ell-1}}+s \leq \nu_j\big\}.$$
Choose $\ell_0$ satisfying
$$\ell_0>\frac{2E(x)}{\epsilon}.$$
We estimate using the gradient flow equation
\begin{eqnarray*}
E(x)&=&\int_{-\infty}^\infty||\partial_s x||^2 ds\\
&\geq&\sum_{\ell=0}^{\ell_0-1}\int_{\nu_{j_\ell}}^{\nu_{j_\ell}+s}||\partial_s x||^2ds
\\
&=&\sum_{\ell=0}^{\ell_0-1}\int_{\nu_{j_\ell}}^{\nu_{j_\ell}+s}\frac{d}{ds} f(x(s))ds
\\
&=&\sum_{\ell=0}^{\ell_0-1}\Big(f(x(\nu_{j_\ell}+s))-
f(x(\nu_{j_\ell}))\Big)\\
&\geq&\frac{\ell_0 \epsilon}{2}\\
&>&E(x).
\end{eqnarray*}
This contradiction shows that the assumption that $x^+$ was nonconstant
had to be wrong. Hence $x^+$ is a critical point and
since $f$ is Morse, the gradient flow line converges at the
positive asymptotic to $x^+$. A completely analogous
reasoning shows that $x$ converges at the negative
asymptotic, too. This proves the Lemma.
\hfill $\square$
\\ \\
Our second assumption is that $(f,g)$ meet the Morse-Smale condition.
We do not suppose that the flow of $\nabla f$ exists
for all times. So instead of assuming that the stable and
unstable manifolds for each pair of critical points of the Morse function
intersect transversely the Morse-Smale condition
has to be rephrased in the assumption that the operator coming from
the linearization of the gradient flow is surjective
as in \cite{schwarz}. In order to recall
this operator we have to introduce some notation.
It follows from Lemma~\ref{limit} that asymptotically each
gradient flow line converges to critical points. For critical
points $x^\pm \in \mathcal{C}$
abbreviate by $\mathcal{H}=\mathcal{H}(x^-,x^+)$ the Hilbert manifold of
$W^{1,2}$-paths from $\mathbb{R}$ to $M$ which converge to $x^\pm$
for $s \to \pm \infty$. Let $\mathcal{E}$ be the bundle over
$\mathcal{H}$ whose fiber at a point $x \in \mathcal{H}$
is given by
$$\mathcal{E}_x=L^2(\mathbb{R},x^* TM).$$
Consider the section
$$\varsigma \colon \mathcal{H} \to \mathcal{E}, \quad
x \mapsto \partial_s x-\nabla f(x).$$
The zero set of this section are gradient flow lines from
$x^-$ to $x^+$. If $x \in \varsigma^{-1}(0)$ there is a canonical
splitting of the tangent space
$$T_x \mathcal{E}=\mathcal{E}_x \otimes T_x \mathcal{H}.$$
Denote by
$$\pi \colon T_x \mathcal{E} \to \mathcal{E}_x$$
the projection along $T_x \mathcal{H}$. The \emph{vertical differential}
at a zero of the section $\varsigma$ at a zero
$x \in \varsigma^{-1}(0)$ is given by
$$D \varsigma(x)= \pi \circ d\varsigma(x) \colon
T_x \mathcal{H}=W^{1,2}(\mathbb{R},x^* TM) \to \mathcal{E}_x.$$
We can now formulate our second hypothesis
\begin{description}
 \item[(H2)] For each $x \in \mathcal{G}$ the operator
$D \varsigma(x)$ is surjective.
\end{description}
\begin{fed}
A tuple $(f,g)$ consisting of a Morse function $f$ and a
Riemannian metric $g$ on the manifold $M$ such that (H1) and (H2) hold is
called a \emph{Morse tuple for $M$}.
\end{fed}
\emph{Remark: }Hypothesis (H1) is actually much more important
than hypothesis (H2) in order to define Morse homology. Even
if transversality fails one can define Morse homology
by using abstract perturbation theory provided compactness
is guaranteed. However, we assume in this paper hypothesis
(H2) so that we can avoid discussions about abstract perturbations.
\\ \\
In the following we assume that we have fixed a Morse
tuple $(f,g)$ on $M$.
Fix further a field  $\mathbb{F}$.
The Morse complex
$$\big(CM_a^b,\partial_a^b\big)=
\big(CM_a^b(f;\mathbb{F}),\partial_a^b(f,g;\mathbb{F})\big)$$
is defined in the following way. The chain group
$$CM_a^b=\mathcal{C}_a^b \otimes \mathbb{F}$$
is the $\mathbb{F}$-vector space
generated by the critical points of $f$ in the action window
$[a,b]$.
Note that $CM_a^b$ is a finite dimensional vector space.
Indeed, it follows from (H1) that the set $\mathcal{C}_a^b$
is compact. Since $f$ is Morse, it is also discrete and hence finite.
The boundary operator $\partial_a^b$ is given by
counting gradient flow lines. For $x^\pm \in \mathcal{C}$ abbreviate
$$\mathcal{G}(x^-,x^+)=\{x \in \mathcal{G}: \lim_{r \to \pm \infty}
r_* x=x^\pm\}$$
If $x^- \neq x^+$ then $\mathbb{R}$ acts freely on
$\mathcal{G}(x^-,x^+)$. Moreover, if the Morse indices satisfy
$\mu(x^-)=\mu(x^+)-1$, then it is well known that it follows
from hypotheses (H1) and (H2) that the quotient
$\mathcal{G}(x^-,x^+)/\mathbb{R}$ is a finite set, see
\cite{schwarz}. In this case
we define the integer
$$m(x^-,x^+)=
\#_\sigma\big(\mathcal{G}(x^-,x^+)/\mathbb{R}\big)$$
where $\#_\sigma$ refers to the signed count of the set.
The sign is
determined by the choice of a coherent orientation for the
moduli spaces of gradient flow lines.
For $c \in \mathcal{C}_a^b$, we put
$$\partial_a^b c=\sum_{\substack{c' \in
\mathcal{C}_a^b\\
\mu(c')=\mu(c)-1}} m(c',c)c'.$$
We define $\partial_a^b$ on $CM_a^b$ by $\mathbb{F}$-linear extension
of the formula above. Again it is well known, see \cite{schwarz}, that
under hypothesis (H1) and (H2) the
homomorphism $\partial_a^b$ is a boundary operator, i.e.
$$\big(\partial_a^b\big)^2=0.$$
Hence we get a graded vector space
$$HM_a^b=HM_a^b(f,g;\Gamma)=
\frac{\mathrm{ker}\partial_a^b}{\mathrm{im}\partial_a^b}.$$
If $a_1 \leq a_2$ we denote by $\underline{\mathcal{C}}^{a_2}_{a_1}$
the set generated by critical points in the half open action
interval $[a_1,a_2)$. In particular, if $a_2$ lies not in
the spectrum of $f$ the graded set $\underline{\mathcal{C}}^{a_2}_{a_1}$
equals $\mathcal{C}^{a_2}_{a_1}$. We abbreviate
\begin{equation}\label{und}
\underline{CM}^{a_2}_{a_1}=\underline{\mathcal{C}}^{a_2}_{a_1}
\otimes \Gamma.
\end{equation}
If $a_1 \leq a_2 \leq b$, then the disjoint union
$$\mathcal{C}_{a_1}^b=\underline{\mathcal{C}}_{a_1}^{a_2} \sqcup
\mathcal{C}_{a_2}^b$$
leads to the direct sum
$$CM_{a_1}^b=\underline{CM}_{a_1}^{a_2} \oplus CM_{a_2}^b.$$
We denote by
$$p^b_{a_2,a_1} \colon CM_{a_1}^b
\to CM_{a_2}^b$$
the projection along $\underline{CM}_{a_1}^{a_2}$. Since the action is
increasing along gradient flow lines the projections commute
with the
boundary operators in the sense that
\begin{equation}\label{proj1}
p^b_{a_2,a_1} \circ \partial_{a_1}^b=
\partial_{a_2}^b \circ p^b_{a_2,a_1}.
\end{equation}
Moreover, for $a_1 \leq a_2 \leq a_3 \leq b$, their composition
obviously meets
\begin{equation}\label{proj2}
p^b_{a_3,a_2} \circ p^b_{a_2,a_1}=p^b_{a_3,a_1}
\end{equation}
and for $a\leq b$
\begin{equation}\label{proj2a}
p^b_{a,a}=\mathrm{id}|_{CM_a^b}.
\end{equation}
It follows from (\ref{proj1}) that
$p^b_{a_2,a_1}$ induces homomorphisms
$$Hp^b_{a_2,a_1} \colon HM_{a_1}^b \to HM_{a_2}^b$$
which satisfy
\begin{equation}\label{proj3}
Hp^b_{a_3,a_2} \circ Hp^b_{a_2,a_1}=Hp^b_{a_3,a_1}.
\end{equation}
by (\ref{proj2}) and
\begin{equation}\label{proj3a}
Hp^b_{a,a}=\mathrm{id}|_{HM_a^b}
\end{equation}
by (\ref{proj2a}).
Similarly, for $a \leq b_1 \leq b_2$ the inclusions
$\mathcal{C}_a^{b_1} \mapsto \mathcal{C}_a^{b_2}$
induce maps
$$i_a^{b_2,b_1} \colon CM_a^{b_1} \to CM_a^{b_2}.$$
Again the inclusions commute with the boundary operators
\begin{equation}\label{in1}
i_a^{b_2,b_1} \circ \partial_a^{b_1}=\partial_a^{b_2}
\circ i_a^{b_2,b_1},
\end{equation}
their composition satisfies
\begin{equation}\label{in2}
i_a^{b_3,b_2} \circ i_a^{b_2,b_1}=i_a^{b_3,b_1}
\end{equation}
for $a \leq b_1 \leq b_2 \leq b_3$ and
\begin{equation}\label{in2a}
i_a^{b,b}=\mathrm{id}|_{CM_a^b}
\end{equation}
for $a \leq b$.
Moreover, inclusions and projections commute in
the sense that if $a_1 \leq a_2 \leq b_1 \leq b_2$, then
\begin{equation}\label{ip1}
i_{a_2}^{b_2,b_1} \circ p^{b_1}_{a_2,a_1}=p^{b_2}_{a_2,a_1}
\circ i_{a_1}^{b_2,b_1}.
\end{equation}
It follows from (\ref{in1}) that $i_a^{b_2,b_1}$ induces
homomorphisms
$$Hi_a^{b_2,b_1} \colon HM_a^{b_1} \to HM_a^{b_2}.$$
By (\ref{in2}) they satisfy
\begin{equation}\label{in3}
Hi_a^{b_3,b_2} \circ Hi_a^{b_2,b_1}=
Hi_a^{b_3,b_1},
\end{equation}
by (\ref{in2a})
\begin{equation}\label{in3a}
Hi_a^{b,b}=\mathrm{id}|_{HM_a^b},
\end{equation}
and by (\ref{ip2})
\begin{equation}\label{ip2}
Hi_{a_2}^{b_2,b_1} \circ Hp^{b_1}_{a_2,a_1}=
Hp^{b_2}_{a_2,a_1} \circ Hi_{a_1}^{b_2,b_1}.
\end{equation}
We can summarize the results of this section in the following
proposition.
\begin{prop}
For a Morse tuple
$(f,g)$ on the manifold $M$ the quadruple
$(CM,p,i,\partial)$ is a bidirect system of chain complexes (see
Section~\ref{ss:bidirect} for the definition).
\end{prop}
\subsection[An axiomatic approach]{An axiomatic approach}\label{axiom}

Following a suggestion of D.\,Salamon we can axiomatize the
results of the previous subsection in the following way.
Via this axiomatized approach one can getTheorem A also in the
infinite dimensional case of Floer homology provided one has the
necessary compactness.

\begin{fed}
A \emph{Floer triple}
$$\mathcal{F}=(\mathcal{C},f,m)$$ consists of a set $\mathcal{C}$,
a function $f \colon \mathcal{C} \to \mathbb{R}$ and a function
$m \colon \mathcal{C} \times \mathcal{C} \to \mathbb{F}$ such that
the following condition holds.
\begin{description}
 \item[(i)] For each $a \leq b$ the set
  $\mathcal{C}^b_a=\{c \in \mathcal{C}: a \leq f(c) \leq b\}$
 is finite.
 \item[(ii)] If $c_1,c_2 \in \mathcal{C}$ and $m(c_1,c_2) \neq 0$,
  then it follows that $f(c_1) < f(c_2)$.
 \item[(iii)] If $c_1,c_3 \in \mathcal{C}$, then
  $\sum_{c_2 \in \mathcal{C}}m(c_1,c_2) m(c_2,c_3)=0$.
\end{description}
\end{fed}
Assertions (i) and (ii) make sure that the sum in
assertion (iii) is finite. Elements of $\mathcal{C}$ are referred
to as critical points, the value of $f$ as their action value and
the number $m(c_1,c_2)$ as the number of gradient flow lines between
$c_1$ and $c_2$. Assertion (i) can then be rephrased by saying that
in each finite action window there are only finitely many critical points,
assertion (ii) says that the action is increasing along gradient flow lines,
and assertion (iii) guarantees that on each action window a boundary
operator can be defined by counting gradient flow lines. As in the
previous subsection one can associate to each Floer triple a
bidirect system of chain complexes. The assumption to have a Morse tuple
in order that Theorem A holds can be generalized to arbitrary Floer
triples.

\section[Algebraic preliminaries]{Algebraic preliminaries}

\subsection[Direct and inverse limits]{Direct and inverse limits}

We first recall that a quasi ordered set is a tuple
$\mathcal{A}=(A,\leq)$ where $A$ is a set and $\leq$ is a
reflexive and transitive binary relation.
More sophisticatedly, one might think of $\mathcal{A}$ as
a category with precisely one morphism from
$a_1$ to $a_2$ whenever $a_1 \leq a_2$. A quasi ordered set
is called partially ordered if the binary
relation is also antisymmetric.

To define direct and inverse limits the notion of
a \emph{direct system} is needed. For the applications we
have in mind we have to work in the category of graded vector spaces.
For simplicity we skip the reference to the grading.
Hence a direct system is a tuple
$$\mathcal{D}=(G,\pi)$$
where $G$ is a family of vector spaces indexed by a
quasi ordered set $\mathcal{A}=(A,\leq)$, i.e.
$$G=\{G_a\}_{a \in A},$$
and
$$\pi=\{\pi_{a_2,a_1}\}_{a_1 \leq a_2;\,a_1,a_2 \in A}$$
is a family
of  homomorphisms
$$\pi_{a_2,a_1} \colon G_{a_1} \to G_{a_2}$$
satisfying
$$\pi_{a,a}=\mathrm{id}|_{G_a},\,\, a \in A, \quad
\pi_{a_3,a_1}=\pi_{a_3,a_2} \circ \pi_{a_2,a_1},\,\,a_1 \leq a_2
\leq a_3.$$
If one thinks of a quasi ordered set as a category, then
a direct system is a functor to the category of vector spaces.
For a direct system the \emph{inverse limit} or just
\emph{limit}
is defined as the vector space
$$\underleftarrow{\lim} G:=\lim_{\substack{\longleftarrow\\
\mathcal{A}}}G
:=\bigg\{\{x_a\}_{a \in A} \in
\prod_{a \in A} G_a: a_1 \leq a_2
\Rightarrow \pi_{a_2,a_1}(x_{a_1})=x_{a_2}\bigg\}.$$
For $a \in A$ let
$$\pi_a \colon \underleftarrow{\lim}G \to G_a$$
be the (not necessarily surjective) projection to the
$a$-th component. These maps  satisfy for $a_1 \leq a_2$
the relation
$$\pi_{a_2}=\pi_{a_2,a_1} \circ \pi_{a_1}.$$
The inverse limit is characterised by
the following universal property. Given a vector space
$H$ and a family of homomorphisms
$\tau_a \colon H \to G_a$ for $a \in A$ which satisfies
$$\tau_{a_2}=\pi_{a_2,a_1}\circ \tau_{a_1}, \quad a_1 \leq a_2,$$
then there exists a unique homomorphism
$\tau \colon H \to \underleftarrow{\lim}G$
such that for any $a \in A$ the following diagram commutes
\begin{equation}
\begin{xy}
 \xymatrix{
   H \ar@{.>}[rr]^{\exists !\,\tau}
  \ar[rd]_{\tau_a}  &   & \underleftarrow{\lim}G
    \ar[dl]^{\pi_a} \\
             & G_a &
}
\end{xy}
\end{equation}
The \emph{direct limit} or \emph{colimit} is constructed
dually to the inverse limit. To make the notation easier
adaptable to our later purposes we denote in the
definition of the direct limit the family of homomorphisms
by $\iota$, i.e.~our direct system reads now
$$\mathcal{D}=(G,\iota).$$
Moreover, the index set is now denoted by
$\mathcal{B}=(B,\leq)$ and subscripts are replaced
by superscripts. For $b \in B$ let
$$\lambda^b \colon G^b \to \bigoplus_{b' \in B} G^{b'}$$
be the $b$-th injection into the sum of the abelian groups
$G^b$. Define the subgroup $S_{\mathcal{D}}$ of
$\bigoplus G^b$ by
$$S_{\mathcal{D}}=
\bigg\{\lambda^{b_2} \iota^{b_2,b_1}(x)-\lambda^{b_1}(x):
x \in G^{b_1},\,\,b_1 \leq b_2\bigg\}.$$
The direct limit is now defined as the vector space
$$\underrightarrow{\lim} G:=\lim_{\substack{\longrightarrow\\
\mathcal{B}}}G
:=\Big(\bigoplus_{b \in B} G^b\Big)\Big/S_{\mathcal{D}}.$$
The direct limit is characterized by the following universal
property dual to the characterization of
the inverse limit. For $b \in B$ let
$$\iota^b \colon G^b \to \underrightarrow{\lim} G$$
be the (not necessarily injective) homomorphism induced from the
inclusion of $G^b$ into the sum $\bigoplus G^b$.
Assume that $H$ is an vector space and
$\tau^b$ for $b \in B$ is a family of homomorphisms
$\tau^b \colon G^b \to H$ satisfying
$$\tau^{b_1}=\tau^{b_2} \circ \iota^{b_2,b_1}, \quad b_1 \leq b_2.$$
Then there exists a unique homomorphism $\tau \colon
\underrightarrow{\lim} G \to H$ such that the following diagram commutes
for any $b \in B$
\[
\begin{xy}
 \xymatrix{
   \underrightarrow{\lim}G \ar@{.>}[rr]^{\exists !\,\tau}
  \ar@{<-}[rd]_{\iota^b}  &   & H
    \ar@{<-}[dl]^{\tau^b} \\
             & G^b &
}
\end{xy}
\]
\subsection[The canonical homomorphism]{The canonical homomorphism}
Direct and inverse limits do not necessarily commute. However, there
is a canonical homomorphism
$$\kappa \colon \underrightarrow{\lim}\underleftarrow{\lim}G
\to \underleftarrow{\lim}\underrightarrow{\lim}G$$
which we describe next. We consider two quasi ordered sets
$\mathcal{A}=(A,\leq)$ and $\mathcal{B}=(B,\leq)$ and
a double indexed family of abelian groups $G_a^b$ with
$a \in A$ and $b \in B$. We suppose that
for every $b \in B$
and every $a_1 \leq a_2 \in A$ there exists a homomorphism
$$\pi^b_{a_2, a_1} \colon G_{a_1}^b \to G_{a_2}^b$$
and for every $a \in A$ and $b_1 \leq b_2 \in B$ there exists a
homomorphism
$$\iota_a^{b_2,b_1} \colon G_a^{b_1} \to G_a^{b_2}$$
such that the following
holds. For any fixed $b \in B$ and any fixed $a \in A$ the tuples
$\big(G^b,\pi^b\big)$ and
$\big(G_a,\iota_a\big)$
are direct systems. Moreover, $\pi$ and $\iota$ are required to commute
in the following sense
$$\iota^{b_2,b_1}_{a_2} \circ \pi^{b_1}_{a_2,a_1}
=\pi^{b_2}_{a_2,a_1} \circ \iota^{b_2,b_1}_{a_1} \colon
G_{a_1}^{b_1} \to G_{a_2}^{b_2},\quad a_1 \leq a_2, \,\,b_1 \leq b_2.$$
This can be rephrased by saying that for every
$a_1 \leq a_2$ and $b_1 \leq b_2$ the square
\begin{equation}\label{pii}
\begin{xy}
 \xymatrix{
  G_{a_1}^{b_1} \ar[r]^{\pi^{b_1}_{a_2,a_1}}
  \ar[d]_{\iota^{b_1,b_2}_{a_1}} & G^{b_1}_{a_2}
  \ar[d]^{\iota_{a_2}^{b_2,b_1}}\\
  G^{b_2}_{a_1} \ar[r]_{\pi_{a_2,a_1}^{b_2}} & G^{b_2}_{a_2}
 }
\end{xy}
\end{equation}
is commutative.
We refer to the triple $\big(G,\pi,\iota\big)$ as a
\emph{bidirect system}.
Due to the commutation relation between $\pi$ and $\iota$
for $a_1 \leq a_2 \in A$ the map
$$\pi_{a_2,a_1} \colon
\lim_{\substack{\longrightarrow\\ \mathcal{B}}}G_{a_1}
\to \lim_{\substack{\longrightarrow\\ \mathcal{B}}}G_{a_2},
\quad \big[\{x^b\}_{b\in B}\big] \mapsto
\big[\{\pi^b_{a_2,a_1}(x^b)\}_{b\in B}\big]$$
is a well defined homomorphism. Analoguously,
for $b_1 \leq b_2 \in B$, we have a well defined homomorphism
$$\iota^{b_2,b_1} \colon \lim_{\substack{\longleftarrow\\
\mathcal{A}}}G^{b_1} \to
\lim_{\substack{\longleftarrow\\ \mathcal{A}}}G^{b_2},
\quad \big\{x_a\big\}_{a \in A} \mapsto
\big\{\iota^{b_2,b_1}_a(x_a)\big\}_{a \in A}.$$
Moreover, both
$(\underrightarrow{\lim}G,\pi)$
and $(\underleftarrow{\lim}G,\iota)$ are
direct systems.
\begin{prop}\label{caho}
For a bidirect system $(G,\pi,\iota)$ there exists for every $b \in B$
a unique homomorphism
$$\kappa^b \colon \lim_{\substack{\longleftarrow\\ \mathcal{A}}}
G^b \to \lim_{\substack{\longleftarrow\\ \mathcal{A}}}
\lim_{\substack{\longrightarrow\\ \mathcal{B}}}G$$
and a unique homomorphism
$$\kappa \colon \lim_{\substack{\longrightarrow\\ \mathcal{B}}}
\lim_{\substack{\longleftarrow\\ \mathcal{A}}}G
\to \lim_{\substack{\longleftarrow\\ \mathcal{A}}}
\lim_{\substack{\longrightarrow\\ \mathcal{B}}}G$$
such that for each $a \in A$ and $b \in B$ the following
diagram commutes
\[
\begin{xy}
 \xymatrix{
G_a^b  \ar[d]^{\iota^b_a}& \ar[l]_{\pi^b_a}
\underleftarrow{\lim}G^b \ar[r]^{\iota^b}
\ar@{.>}[d]^{\exists!\,\kappa^b}& \underrightarrow{\lim}
\underleftarrow{\lim}G
\ar@{.>}[ld]^/-.8em/{\exists!\, \kappa}\\
\underrightarrow{\lim} G_a  &  \ar[l]_{\pi_a}
\underleftarrow{\lim}\underrightarrow{\lim}G
}
\end{xy}
\]
\end{prop}
\textbf{Proof: }A straightforward computation shows that
for $a_1 \leq a_2 \in A$ and $b \in B$ the formula
$$\iota_{a_2}^b\pi_{a_2}^b=\pi_{a_2,a_1}
\iota_{a_1}^b\pi_{a_1}^b$$
holds. Hence existence and uniqueness of the $\kappa^b$
for $b \in B$ follows from the universal property of
the inverse limit. For $b_1 \leq b_2 \in B$ and $a \in A$
one computes using the already establishes commutativity in
the left square that
$$\pi_a\kappa^{b_2}\iota^{b_2,b_1}=\iota_a^{b_1}\pi_a^{b_1}=
\pi_a \kappa^{b_1}.$$
Using uniqueness we conclude that
$$\kappa^{b_2}\iota^{b_2,b_1}=\kappa^{b_1}.$$
Now existence and uniqueness of $\kappa$ follows from the universal
property of the direct limit. \hfill $\square$
\\ \\
Conditions under which the canonical homomorphism $\kappa$
is an isomorphism were obtained by B.\,Eckmann and
P.\,Hilton in \cite{eckmann-hilton} and by
A.\,Frei and J.\,Macdonald in \cite{frei-macdonald}. We first
remark that in general $\kappa$ is neither necessarily
injective nor surjective. The example at the end of
section~\ref{intro} shows that injectivity might fail. An
example were surjectivity fails is described in
\cite[p.\,117]{eckmann-hilton}.

To describe the result
of A.\,Frei and J.\,Macdonald we need the following terminology.
To a square
\begin{equation}\label{ehs}
\begin{xy}
 \xymatrix{
  A \ar[r]^{\phi_{BA}} \ar[d]_{\phi_{CA}} & B \ar[d]^{\phi_{DB}}\\
  C \ar[r]_{\phi_{DC}} & D
 }
\end{xy}
\end{equation}
we can associate the sequence
\begin{equation}\label{eh}
A \stackrel{\{\phi_{BA},\phi_{CA}\}}{\longrightarrow}
B \oplus C \stackrel{\langle \phi_{DB},-\phi_{DC}\rangle}
{\longrightarrow} D.
\end{equation}
The square (\ref{ehs}) is commutative precisely if
the sequence (\ref{eh}) is a complex. A commutative square
is now called \emph{exact}, \emph{cartesian}, \emph{cocartesian},
or \emph{bicartesian} iff the corresponding sequence is
exact, left exact, right exact, or a short exact sequence.
\\
We further recall that a quasi ordered set $\mathcal{A}=(A,\leq)$
is \emph{upward directed} if for any $a, a' \in A$ there exists
$a'' \in A$ such that $a \leq a''$ and $a' \leq a''$. Dually
it is called \emph{downward directed} if for any
$a, a' \in A$ there exists $a'' \in A$ such that
$a'' \leq a$ and $a'' \leq a'$.
\\
The following Theorem follows from
\cite[Theorem 5.6]{frei-macdonald}.
\begin{thm}\label{armin}
Assume that $\mathcal{A}$ is upward directed, $\mathcal{B}$
is downward directed and the commutative square (\ref{pii})
is cartesian for any $a_1 \leq a_2$ and $b_1 \leq b_2$.
Then the canonical homomorphism $\kappa
\colon \underrightarrow{\lim} \underleftarrow{\lim}G
\to \underleftarrow{\lim}\underrightarrow{\lim}G$
is an isomorphism.
\end{thm}

\subsection[Bidirect systems of chain complexes]{Bidirect systems
of chain complexes}\label{ss:bidirect}

A bidirect system of chain
complexes is a quadruple
$$\mathcal{Q}=(C,p,i,\partial)$$
where $(C,p,i)$ is a bidirect system which in addition is endowed for
each $a \in A$ and $b \in B$ with a boundary operator
$$\partial^b_a \colon C^b_a \to C^b_a$$
which commutes with $i$ and $p$ in the sense of (\ref{proj1}) and
(\ref{in1}). If
$$HC^b_a=\frac{\mathrm{ker}\partial^b_a}{\mathrm{im}\partial^b_a}$$
are the homology groups, and
$Hp^{b}_{a_2,a_1}$ and $Hi^{b_2,b_1}_a$ are the induced maps
on homology the triple $(HC,Hp,Hi)$ is a bidirect system. As
in the previous subsection we
let
$$\kappa \colon \underrightarrow{\lim}\underleftarrow{\lim}HC
\to \underleftarrow{\lim}\underrightarrow{\lim}HC$$
be the canonical homomorphism on homology level. We refer to
$$k \colon \underrightarrow{\lim}\underleftarrow{\lim}C
\to \underleftarrow{\lim}\underrightarrow{\lim}C$$
as the canonical homomorphism on chain level. Since
$\partial$ commutes with $i$ and $p$ we obtain an induced map
$$Hk \colon H\big(\underrightarrow{\lim}\underleftarrow{\lim}C\big)
\to H\big(\underleftarrow{\lim}\underrightarrow{\lim}C\big).$$
Moreover, for $a \in A$ and $b \in B$
the maps
$$Hi^b_a \colon HC^b_a \to H(\underrightarrow{\lim}C_a)$$
satisfy for $b_1 \leq b_2$
$$Hi^{b_1}_a=Hi^{b_2}_a \circ Hi^{b_2,b_1}_a$$
and hence by the universal property of the direct limit there
exists a unique map
$$\mu_a \colon \underrightarrow{\lim}HC_a \to
H(\underrightarrow{\lim}C_a)$$
such that for any $b \in B$ the diagram
\[
\begin{xy}
 \xymatrix{
   \underrightarrow{\lim}HC_a \ar[rr]^{\mu_a}
  \ar@{<-}[rd]_{\iota^b_a}  &   & H(\underrightarrow{\lim}C_a)
    \ar@{<-}[dl]^{Hi_a^b} \\
             & HC_a^b &
}
\end{xy}
\]
commutes. Taking inverse limits of this diagram and using
functoriality of the inverse limit gives a commutative diagram
\[
\begin{xy}
 \xymatrix{
   \underleftarrow{\lim}\underrightarrow{\lim}HC
   \ar[rr]^{\underleftarrow{\lim}\mu}
  \ar@{<-}[rd]_{\underleftarrow{\lim}\iota^b}  &   &
  \underleftarrow{\lim}H(\underrightarrow{\lim}C)
    \ar@{<-}[dl]^{\underleftarrow{\lim}Hi^b} \\
             & \underleftarrow{\lim}HC^b &
}
\end{xy}
\]
Taking first inverse limits of the chain complexes and
applying the procedure above gives a map
$$\mu \colon \underrightarrow{\lim}H(\underleftarrow{\lim}C)
\to H(\underrightarrow{\lim}\underleftarrow{\lim}C)$$
which is uniquely characterised by the commutativity of the
following diagram
\[
\begin{xy}
 \xymatrix{
   \underrightarrow{\lim}H(\underleftarrow{\lim}C) \ar[rr]^{\mu}
  \ar@{<-}[rd]_{\iota^b}  &   & H(\underrightarrow{\lim}
   \underleftarrow{\lim} C)
    \ar@{<-}[dl]^{Hi^b} \\
             & H(\underleftarrow{\lim}C^b) &
}
\end{xy}
\]
Similarly by using the universal property of the inverse limit
we obtain for each $b \in B$ a map
$$\nu^b \colon H(\underleftarrow{\lim}C^b) \to \underleftarrow{\lim}
HC^b$$
such that for each $a \in A$ the diagram
\begin{equation}\label{nu}
\begin{xy}
 \xymatrix{
   H(\underleftarrow{\lim}C^b) \ar[rr]^{\nu^b}
  \ar[rd]_{Hp^b_a}  &   & \underleftarrow{\lim}HC^b
    \ar[dl]^{\pi_a^b} \\
             & HC_a^b &
}
\end{xy}
\end{equation}
commutes. Taking the direct limit of this diagram we get the commutative
diagram
\[
\begin{xy}
 \xymatrix{
   \underrightarrow{\lim}H(\underleftarrow{\lim}C)
    \ar[rr]^{\underrightarrow{\lim}\nu}
  \ar[rd]_{\underrightarrow{\lim}Hp^b}  &   &
\underrightarrow{\lim}\underleftarrow{\lim}HC
    \ar[dl]^{\underrightarrow{\lim}\pi^b} \\
             & \underrightarrow{\lim}HC^b &
}
\end{xy}
\]
Applying the direct limit already on chain level we obtain a map
$$\nu \colon H(\underleftarrow{\lim}\underrightarrow{\lim}C)
\to \underleftarrow{\lim}H(\underrightarrow{\lim}C)$$
such that the following diagram commutes
\[
\begin{xy}
 \xymatrix{
   H(\underleftarrow{\lim}\underrightarrow{\lim}C) \ar[rr]^{\nu}
  \ar[rd]_{Hp_a}  &   & \underleftarrow{\lim}H(\underrightarrow{\lim}C)
    \ar[dl]^{\pi_a} \\
             & H(\underrightarrow{\lim}C_a) &
}
\end{xy}
\]
We summarize the plethora of maps we passed by in the diagram
\[
\begin{xy}
 \xymatrix{
H(\underrightarrow{\lim}\underleftarrow{\lim}C)
\ar[d]^{Hk}& \ar[l]_{\mu}
\underrightarrow{\lim}H(\underleftarrow{\lim}C)
\ar[r]^{\underrightarrow{\lim}\nu}
& \underrightarrow{\lim}
\underleftarrow{\lim}HC
\ar[d]^{\kappa}\\
H(\underleftarrow{\lim}\underrightarrow{\lim}C)   \ar[r]^{\nu} &
\underleftarrow{\lim}H(\underrightarrow{\lim}C) &
\ar[l]_{\underleftarrow{\lim}\mu}
\underleftarrow{\lim}\underrightarrow{\lim}HC
}
\end{xy}
\]
We do not know if the diagram above always commutes. But we
make now an assumption on the bidirect system of chain complexes
which guarantees commutativity of the diagram above.
\begin{fed}\label{tame}
A bidirect system of chain complexes is called \emph{tame}
if for any $a \in A$ and any $b \in B$ the maps $\mu_a$ and
$\nu^b$ as well as the map $\mu$ are isomorphisms.
\end{fed}
We point out that for a tame bidirect system of chain complexes
the map $\nu$ does not need to be an isomorphism. However,
by functoriality of the inverse and direct limits the maps
$\underleftarrow{\lim}\mu$ and $\underrightarrow{\lim}\nu$ are
isomorphisms, too. Hence we can define maps
$$\rho \colon H(\underleftarrow{\lim}\underrightarrow{\lim}C)
\to \underleftarrow{\lim}\underrightarrow{\lim}HC, \quad
\rho=(\underrightarrow{\lim}\nu) \circ \mu^{-1}$$
and
$$\sigma \colon H(\underleftarrow{\lim}\underrightarrow{\lim}C)
\to \underleftarrow{\lim}\underrightarrow{\lim}HC, \quad
\sigma=(\underleftarrow{\lim}\mu)^{-1} \circ \nu.$$
In particular, the previous diagram simplifies to
\begin{equation}\label{maindia}
\begin{xy}
 \xymatrix{
H(\underrightarrow{\lim}\underleftarrow{\lim}C)
\ar[d]^{Hk}
\ar[r]^{\rho}
& \underrightarrow{\lim}
\underleftarrow{\lim}HC
\ar[d]^{\kappa}\\
H(\underleftarrow{\lim}\underrightarrow{\lim}C)   \ar[r]^{\sigma} &
\underleftarrow{\lim}\underrightarrow{\lim}HC
}
\end{xy}
\end{equation}
\begin{prop}\label{ta}
Assume that the bidirect system is tame. Then the
diagram (\ref{maindia})
commutes and $\rho$ is an isomorphism.
\end{prop}
\textbf{Proof: }That $\rho$ is an isomorphism is clear since
as we observed above $\underrightarrow{\lim}\nu$ is an isomorphism.
We show commutativity in two steps.
\\ \\
\textbf{Step~1: }\emph{For every $b \in B$ the following diagram
commutes}
\[
\begin{xy}
 \xymatrix{
H(\underleftarrow{\lim}C^b)
\ar[d]^{Hk^b}
\ar[r]^{\nu^b}
& \underleftarrow{\lim}HC^b
\ar[d]^{\kappa^b}\\
H(\underleftarrow{\lim}\underrightarrow{\lim}C)   \ar[r]^{\sigma} &
\underleftarrow{\lim}\underrightarrow{\lim}HC
}
\end{xy}
\]
To prove Step~1 we enlarge the diagram to the following one
\[
\begin{xy}
 \xymatrix{
H(\underleftarrow{\lim}C^b)
\ar[d]^{Hk^b}
\ar[r]^{\nu^b}
& \underleftarrow{\lim}HC^b
\ar[d]^{\kappa^b} \ar[r]^{\pi^b_a} & HC^b_a \ar[d]^{\iota^b_a}
\ar[dr]^/1.2em/{Hi^b_a}\\
H(\underleftarrow{\lim}\underrightarrow{\lim}C)   \ar[r]^{\sigma} &
\underleftarrow{\lim}\underrightarrow{\lim}HC \ar[r]^{\pi_a}
& \underrightarrow{\lim}HC_a \ar[r]^{\mu_a} & H(\underrightarrow{\lim}C_a)
}
\end{xy}
\]
The triangle on the right and the middle square commute. We claim
that the exterior square also commutes. Indeed, this square is
obtained by applying the homology functor to the commutative square
\[
\begin{xy}
 \xymatrix{
\underleftarrow{\lim}C^b
\ar[d]^{k^b}
\ar[r]^{p^b_a}
& C^b_a \ar[d]^{i^b_a}\\
\underleftarrow{\lim}\underrightarrow{\lim}C   \ar[r]^{p_a} &
\underrightarrow{\lim}C_a
}
\end{xy}
\]
Using the fact that $\nu^b$ and $\mu_a$ are isomorphisms we conclude that
the diagram
\[
\begin{xy}
 \xymatrix{
\underleftarrow{\lim}HC^b
\ar@<2pt>[d]^{\kappa^b}
\ar@<-2pt>[d]_{\sigma \circ Hk^b \circ(\nu^b)^{-1}}
\ar[r]^{\pi^b_a}
& C^b_a \ar[d]^{\iota^b_a}\\
\underleftarrow{\lim}\underrightarrow{\lim}HC   \ar[r]^{\pi_a} &
\underrightarrow{\lim}HC_a
}
\end{xy}
\]
is commutative for both arrows. But by Proposition~\ref{caho}
the map $\kappa^b$ is unique with this property. Hence
$$\kappa^b=\sigma \circ HK^b \circ (\nu^b)^{-1}$$
and Step~1 follows.
\\ \\
\textbf{Step~2: } \emph{The diagram (\ref{maindia}) commutes.}
\\ \\
For $b \in B$ we enlarge diagram (\ref{maindia}) to the diagram
\[
\begin{xy}
 \xymatrix{
H(\underrightarrow{\lim}\underleftarrow{\lim}C)
\ar[d]^{Hk}
& \ar[l]_{\rho^{-1}} \underrightarrow{\lim}
\underleftarrow{\lim}HC
\ar[d]^{\kappa} &  \ar[l]_{\iota^b} \underleftarrow{\lim}HC^b
\ar[ld]^/-1.2em/{\kappa^b}\\
H(\underleftarrow{\lim}\underrightarrow{\lim}C)   \ar[r]^{\sigma} &
\underleftarrow{\lim}\underrightarrow{\lim}HC
& H(\underleftarrow{\lim}C^b) \ar@/^/[ll]^{Hk^b} \ar[u]_{\nu^b}
}
\end{xy}
\]
The exterior square is obtained by applying the homology functor
to the commutative triangle
\[
\begin{xy}
 \xymatrix{
   \underrightarrow{\lim}\underleftarrow{\lim}C \ar@{<-}[rr]^{\iota^b}
  \ar[rd]_{k}  &   & \underleftarrow{\lim}C^b
    \ar[dl]^{k^b} \\
             & \underleftarrow{\lim}\underrightarrow{\lim}C &
}
\end{xy}
\]
and is therefore commutative. Hence using Step~1 and the assumption
that $\nu^b$ is an isomorphism we deduce that the
diagram
\[
\begin{xy}
 \xymatrix{
\underrightarrow{\lim}\underleftarrow{\lim}HC
\ar@<2pt>[d]^{\kappa}
\ar@<-2pt>[d]_{\sigma \circ Hk \circ \rho^{-1}}
&  \ar[l]_{\iota^b} \underleftarrow{\lim}HC^b
\ar[ld]^/-1.2em/{\kappa^b}\\
\underleftarrow{\lim}\underrightarrow{\lim}HC
}
\end{xy}
\]
is commutative for both arrows. Again by Proposition~\ref{caho}
we conclude that
$$\kappa=\sigma \circ Hk \circ \rho^{-1}.$$
This finishes the proof of Step~2 and hence of the proposition.
\hfill $\square$

\subsection[The Mittag-Leffler condition]{The Mittag-Leffler
condition}

Given a direct system of chain complexes $(C,p,\partial)$
there is a canonical map $\nu \colon H(\underleftarrow{\lim}C)
\to \underleftarrow{\lim}HC$ defined as in (\ref{nu}).
An important tool to study surjectivity and bijectivity
properties of the map $\nu$ is the Mittag-Leffler condition.
Following A.\,Grothendieck, see \cite[(13.1.2)]{grothendieck},
this condition reads as follows.

\begin{fed}A direct system $(G,\pi)$ of
vector spaces indexed on the quasi-ordered set
$(\mathbb{R},\leq)$ is said to satisfy
the \emph{Mittag-Leffler condition} if for any $a \in A$ there
exists $a'=a'(a) \leq a$ such that for any $a'' \leq a'$ the following
holds
$$\mathrm{im} \pi_{a,a''}=\mathrm{im}\pi_{a,a'} \subset G_a.$$
\end{fed}
The following lemma gives two criteria under which the Mittag-Leffler
condition holds true.
\begin{lemma}\label{ml}
The Mittag-Leffler condition holds in the following two cases.
\begin{description}
\item[(i)] For every $a_1 \leq a_2$ the homomorphism $\pi_{a_2,a_1}
\colon G_{a_1} \to G_{a_2}$ is surjective.
\item[(ii)] For any $a \in \mathbb{R}$ the vector space
$G_a$ is finite dimensional.
\end{description}
\end{lemma}
\textbf{Proof: } That the Mittag-Leffler condition holds in
case (i) is obvious. To show that it holds in case (ii) we first
observe that the relation $\pi_{a,a''}=\pi_{a,a'} \circ \pi_{a',a''}$
for $a'' \leq a' \leq a$ implies that
\begin{equation}\label{increas}
\mathrm{im} \pi_{a,a''} \subset \mathrm{im} \pi_{a,a'} \subset G_a.
\end{equation}
Using that $G_a$ is finite dimensional the function
$$\varrho_a \colon (-\infty,a] \to \mathbb{N} \cup \{0\}, \quad
a' \mapsto \mathrm{dim}\big(\mathrm{im}\pi_{a,a'}\big)$$
is well-defined and it is monotone increasing by (\ref{increas}).
Since it is bounded from below and takes only discrete values there
exists
$$m_a=\mathrm{min} \varrho_a \in \mathbb{N} \cup \{0\}.$$
We choose $a'=a'(a)$ in such a way that
$$\varrho_a(a')=m_a.$$
With this choice it follows that for every $a'' \leq a'$
it holds that
$$\mathrm{dim} \big(\mathrm{im}\pi_{a,a''}\big)=\mathrm{dim}
\big(\mathrm{im}\pi_{a,a'}\big).$$
Hence by (\ref{increas}) we get
$$\mathrm{im}\pi_{a'',a}=\mathrm{im}\pi_{a',a}$$
which finishes the proof of the Mittag-Leffler condition. \hfill $\square$
\\ \\
For the following theorem, see
\cite[Proposition~13.2.3]{grothendieck} or
\cite[Proposition~3.5.7, Theorem~3.5.8]{weibel}.
\begin{thm}\label{tnu}
Assume that $(C,p,\partial)$ is a direct system of chain complexes
indexed on the set $(\mathbb{R},\leq)$.
If $(C,p)$ satisfies the Mittag-Leffler condition, then the
homomorphism $\nu \colon H(\underleftarrow{\lim}C) \to
\underleftarrow{\lim}HC$ is surjective. If in addition
$(HC,Hp)$ satisfies the Mittag-Leffler condition, too, then
$\nu$ is an isomorphism.
\end{thm}
Under the assumptions of Theorem~\ref{tnu}
if $(C,p)$ satisfies the Mittag-Leffler condition, then the
kernel of $\nu$ can be described with the help of the
first derived functor $\underleftarrow{\lim}^1$
of the inverse limit. If $(G,\pi)$ is a direct system
of abelian groups indexed on the real line,
$\underleftarrow{\lim}^1 G$ can be described in the following
way. Choose a sequence $a_j \in \mathbb{R}$ such that
$a_{j+1} \leq a_j$ for every $j \in \mathbb{N}$ and
$a_j$ converges to $-\infty$, i.e. $\{a_j\}_{j \in \mathbb{N}}$
is a cofinal
sequence in $\mathbb{R}$. Consider the map
$$\Delta \colon \prod_{j=1}^\infty G_{a_j}
\to \prod_{j=1}^\infty G_{a_j}, \quad
\big\{x_{a_j}\big\}_{j \in \mathbb{N}} \mapsto
\big\{x_{a_j}-\pi_{a_j,a_{j+1}}
(x_{a_{j+1}})\big\}_{j \in \mathbb{N}}$$
and set
$$\underleftarrow{\lim}^1G=\mathrm{coker} \Delta.$$
It is straightforward to check that $\underleftarrow{\lim}^1 G$
only depends on the choice of the cofinal sequence up to
canonical isomorphism. For a graded abelian group $G$ and
$n \in \mathbb{Z}$ let
$G[n]$ be the graded group obtained from $G$ by shifting the
grading by $n$. Theorem~\ref{tnu} follows from
the following exact sequence
\begin{equation}\label{milnor}
0 \to \underleftarrow{\lim}^1 HC[1]\rightarrow
H(\underleftarrow{\lim}C) \stackrel{\nu}{\rightarrow}
 \underleftarrow{\lim}HC
\rightarrow 0
\end{equation}
and the fact that the Mittag-Leffler condition implies
the vanishing of $\underleftarrow{\lim}^1$. The sequence
(\ref{milnor}) is also known as \emph{Milnor sequence}
since it appeared in a slightly different context in
the work of Milnor, see \cite{milnor}.
\\ \\
\emph{Remark: } One can also define higher derived functors
$\underleftarrow{\lim}^n$ of the inverse limit. This was
carried out by J.\,Roos in \cite{roos} and
G.\,N\"obeling in \cite{noebeling}. However, if the direct system
is indexed on the reals the functors
$\underleftarrow{\lim}^n$ vanish for $n \geq 2$.

\section[Proof of Theorem A]{Proof of Theorem A}

Let $(CM,p,i,\partial)$ be the bidirect system of
chain complexes associated to a Morse tuple $(f,g)$
on a manifold $M$ or more generally to a Floer triple
$\mathcal{F}=(\mathcal{C},f,m)$.
Recall from Definition~\ref{tame} the notion of a tame
bidirect system of chain complexes. We need the following
Lemma.
\begin{lemma}\label{ta2}
The bidirect system $(CM,p,i,\partial)$ is tame.
\end{lemma}
\textbf{Proof: }Since $\mathbb{R}$ is upward directed
the direct limit functor commutes with the homology functor
\cite[Theorem IV.7]{spanier}. Consequently the homomorphism
$\mu$ and the homomorphisms $\mu_a$ for any $a \in \mathbb{R}$
are isomorphisms. Because the projections $p^b_{a_2,a_1}$
are surjective it follows from assertion (i) in
Lemma~\ref{ml} that for any $b \in \mathbb{R}$
the direct system of abelian groups $(CM^b,p^b)$ satisfies
the Mittag-Leffler condition. Since all the vector spaces
$HM^b_a$ are finite dimensional assertion (ii) of Lemma~\ref{ml}
implies that the
direct system $(HM^b,Hp^b)$ satisfies the Mittag-Leffler
condition, too. Hence it follows from Theorem~\ref{tnu}
that the homomorphisms $\nu^b$ for any $b \in \mathbb{R}$ are
also isomorphisms. This proves that
$(CM,p,i,\partial)$ is tame. \hfill $\square$
\\ \\
In view of Proposition~\ref{ta} and Lemma~\ref{ta2} the
following diagram commutes and $\rho$ is an isomorphism
\begin{equation}\label{maindia2}
\begin{xy}
 \xymatrix{
H(\underrightarrow{\lim}\underleftarrow{\lim}CM)
\ar[d]^{Hk}
\ar[r]^{\rho}
& \underrightarrow{\lim}
\underleftarrow{\lim}HM=\overline{HM}
\ar@<-12pt>[d]^{\kappa}\\
H(\underleftarrow{\lim}\underrightarrow{\lim}CM)   \ar[r]^{\sigma} &
\underleftarrow{\lim}\underrightarrow{\lim}HM=\underline{HM}
}
\end{xy}
\end{equation}
For the following Lemma recall that $HM$ is the Morse homology
obtained by taking the Novikov completion of the chain groups
$CM^b_a$.
\begin{lemma}\label{isom}
The homomorphism $k$ and $Hk$ are isomorphisms and
\begin{equation}\label{novi}
H(\underleftarrow{\lim}\underrightarrow{\lim}CM)=
H(\underrightarrow{\lim}\underleftarrow{\lim}CM)=HM.
\end{equation}
\end{lemma}
\textbf{Proof: }If $k$ is an isomorphism, then $Hk$
obviously is an isomorphism, too. To see that
$k$ is an isomorphism observe that the elements of both
$\underleftarrow{\lim}\underrightarrow{\lim}CM$ and
$\underrightarrow{\lim}\underleftarrow{\lim}CM$
are given by Novikov sums
$$\xi=\sum_{c \in \mathcal{C}} \gamma_c c, \quad
\gamma_c \in \mathbb{F}, \quad
\#\{c \in \mathcal{C}: \gamma_c \neq 0,\,\, f(c)>b\}<\infty
,\,\,\forall\,\,b \in \mathbb{R}.$$
This additionally implies the second equality
in (\ref{novi}). \hfill $\square$
\\ \\
Before continuing with the proof of Theorem~A we remark
that the fact that $k$ is an isomorphism can also be
deduced from Theorem~\ref{armin} in view of the following
Lemma.
\begin{lemma}
For the bidirect system $(CM,p,i,\partial)$ each diagram
(\ref{pii}) is bicartesian and hence in particular cartesian.
\end{lemma}
\textbf{Proof: }We have to show that for each $a_1 \leq a_2$
and $b_1 \leq b_2$ the sequence
\begin{equation}
CM^{b_1}_{a_1}
\stackrel{\big\{i_{a_1}^{b_2,b_1},p_{a_2,a_1}^{b_1}\big\}}
{\longrightarrow}
CM^{b_2}_{a_1} \oplus CM^{b_1}_{a_2}
\stackrel{\big\langle p_{a_2,a_1}^{b_2},-i_{a_2}^{b_2,b_1}
\big\rangle}
{\longrightarrow} CM^{b_2}_{a_2}
\end{equation}
is short exact. Since $i_{a_1}^{b_2,b_1}$ is injective
the first map is an injection and since
$p_{a_2,a_1}^{b_2}$ is surjective the second map is a
surjection. It remains to show exactness. Let
$$\Delta_{a_2}^{b_1} \subset CM_{a_2}^{b_1} \oplus
CM_{a_2}^{b_1}$$
be the diagonal. Via the embedding
$CM_{a_2}^{b_1} \oplus CM_{a_2}^{b_1} \hookrightarrow
CM_{a_1}^{b_2} \oplus CM_{a_2}^{b_1}$ we think of
$\Delta_{a_2}^{b_1}$ as a subvectorspace of
$CM_{a_1}^{b_2} \oplus CM_{a_2}^{b_1}$. Recall
the notation $\underline{CM}_{a_1}^{a_2}$ from (\ref{und}).
We then have
$$\mathrm{im}\big\{i_{a_1}^{b_2,b_1},p_{a_2,a_1}^{b_1}\big\}
=\Delta_{a_2}^{b_1} \cup
\big(\underline{CM}_{a_1}^{a_2} \oplus \{0\}\big)
=\mathrm{ker}\big\langle p_{a_2,a_1}^{b_2},-i_{a_2}^{b_2,b_1}
\big\rangle.$$
This shows exactness and hence the lemma is proved. \hfill $\square$
\\ \\
\textbf{End of proof of Theorem~A: } Setting
$\overline{\rho}=\rho$ and $\underline{\rho}=\sigma \circ Hk$
we conclude from the diagram (\ref{maindia2}) using Lemma~\ref{isom}
that the diagram (\ref{dia}) is commutative with
$\overline{\rho}$ an isomorphism. It remains to show
that $\underline{\rho}$ is surjective. Using the formula
$$\underline{\rho}=\sigma \circ Hk=(\underleftarrow{\lim}\mu)^{-1}
\circ \nu \circ Hk$$
and the fact that $\underleftarrow{\lim}\mu$ and $Hk$ are
isomorphisms we are reduced to show that $\nu$ is surjective.
Since for any $a_1 \leq a_2$ the homomorphism
$\underrightarrow{\lim}p_{a_2,a_1}$ is surjective
we conclude that the bidirect system
$(\underrightarrow{\lim}CM,\underrightarrow{\lim}p)$
satisfies the Mittag-Leffler condition. Hence
it follows again from Theorem~\ref{tnu} that
$\nu$ is surjective. We are done with the proof of
Theorem~A. \hfill $\square$
\\ \\
\emph{Remark: } Using Milnor's exact sequence (\ref{milnor})
one observes that the kernel of the canonical homomorphism
$\kappa$ is given by
$$\mathrm{ker}\kappa=\underleftarrow{\lim}^1
H(\underrightarrow{\lim}CM)=\underleftarrow{\lim}^1
\underrightarrow{\lim}HM.$$

\appendix

\section[Integer coefficients]{Integer coefficients}

The homomorphism $\overline{\rho} \colon HM \to \overline{HM}$
need not be an isomorphism any more if one
uses integer coefficients. We show this in an example.
We consider the following Floer triple
$\mathcal{F}=(\mathcal{C},f,m)$.
The critical set $\mathcal{C}$ is given by
$$\mathcal{C}=\big\{\overline{c}_n: n \in \mathbb{N} \cup\{0\}
\big\} \cup \big\{\underline{c}_n: n \in \mathbb{N}\big\}.$$
The function $f$ satisfies
$$f(\overline{c}_n)=-n, \quad f(\underline{c}_n)=-n-1$$
and the nonvanishing entries of $m$ are
$$m(\underline{c}_n,\overline{c}_{n-1})=1, \quad
m(\underline{c}_n,\overline{c}_n)=-2, \qquad
n \in \mathbb{N}.$$
We point out again that in the following theorem we use
integer coefficients.
\begin{thm}
For the Floer triple $\mathcal{F}$ as above,
$\overline{HM}=0$, but $HM \neq 0$.
\end{thm}
\textbf{Proof: }We prove the theorem in three steps. For
$n \in \mathbb{N} \cup \{0\}$ we use the abbreviation
$$\gamma_n=\sum_{j=0}^{n}2^{n-j}\overline{c}_j.$$
\textbf{Step\,1: } \emph{For $b \geq 0$ and $a\leq -1$ with
$k=\lfloor -a \rfloor$ we have}
$$HM_a^b=\mathbb{Z}[\gamma_{k-1}]
\oplus \mathbb{Z}[\gamma_k].$$
We first observe that the chain group is given by
$$CM^b_a=\bigoplus_{j=0}^k \mathbb{Z}\overline{c}_j \oplus
\bigoplus_{j=1}^{k-1} \mathbb{Z}\underline{c}_j.$$
We claim that
\begin{equation}\label{geg1}
\mathrm{im}\partial^b_a=\bigoplus_{j=1}^{k-1}
\mathbb{Z}\underline{c}_j.
\end{equation}
It is clear that the left-hand side is contained in the right-hand
side since there are no gradient flow lines starting from
a critical point $\overline{c}_n$. To see the other inclusion,
observe that
\begin{equation}\label{geg2}
\underline{c}_n=\partial^b_a\gamma_{n-1},
\quad n \in \{1,\ldots,k-1\}
\end{equation}
which implies (\ref{geg1}).
We next show that
\begin{equation}\label{geg3}
\mathrm{ker}\partial^b_a=\mathbb{Z} \gamma_{k-1}
\oplus \mathbb{Z} \gamma_k \oplus \mathrm{im}\partial^b_a.
\end{equation}
It is straightforward to check that the righthand side is
contained in the kernel of the boundary operator. To see
the other inclusion we observe that
$\{\gamma_0,\ldots, \gamma_k\}$ is another
$\mathbb{Z}$-basis of the free abelian group
$\bigoplus_{j=0}^k \mathbb{Z}\overline{c}_j$. Indeed,
the two bases are related by an upper triangular matrix
with diagonal entries one. In particular, the determinant
of this matrix is one. Assertion (\ref{geg3})
therefore follows from (\ref{geg2}). Step\,1 is an
immediate consequence of (\ref{geg3}).
\\ \\
\textbf{Step\,2: } \emph{$\overline{HM}=0$.}
\\ \\
We first prove that for $b \geq 0$ we have
\begin{equation}\label{feh2}
\lim_{\substack{\longleftarrow\\
a \to -\infty}}HM_a^b=0.
\end{equation}
To see that assume that
$$x \in \lim_{\substack{\longleftarrow\\
a \to -\infty}}HM_a^b.$$ Then
$$x=\{x_a\}_{a \leq b}$$
where $x_a \in HM_a^b$ and for $a_1 \leq a_2 \leq b$ the equation
$$Hp^b_{a_2,a_1}(x_{a_1})=x_{a_2}$$ holds.
We have to show that
$$
x_a=0, \quad a \leq b.
$$
This is clear if $a>0$ since in this case $HM_a^b=0$, because there
are no critical points of positive action. If $a \in (-1,0]$ then
$$HM_a^b=\mathbb{Z}[\overline{c}_0]=\mathbb{Z}[\gamma_0].$$
and hence there exists $n_a \in \mathbb{Z}$ such that
\begin{equation}\label{ach1}
   x_a=n_a [\gamma_0], \quad a \in (-1,0].
\end{equation}
Since there are no critical points in the action window $(-1,0)$ we
conclude
$$n_a=n_0, \quad a \in (-1,0]$$
If $a \leq -1$ then by Step\,1 there exist $n^1_a, n^2_a \in
\mathbb{Z}$ such that
\begin{equation}\label{ach2}
 x_a=n^1_a [\gamma_{\lfloor -a
\rfloor-1}]+n^2_a [\gamma_{\lfloor -a \rfloor}], \quad a \leq -1.
\end{equation}
Again since for each $k \in \mathbb{N}$ there are no critical points
in the action window $(-k-1,-k)$ we conclude that
$$n^1_a=n^1_{\lceil a \rceil},\,\,n^2_a=n^2_{\lceil a \rceil}, \quad
a \leq -1.$$ Hence to prove (\ref{feh2}) we are left with showing
\begin{equation}\label{ups}
n_0=0, \qquad n^1_{-k}=n^2_{-k}=0, \quad k \in \mathbb{N}.
\end{equation}
For $k \in \mathbb{N} \cup \{0\}$ and $\ell \in \mathbb{N}$ we
compute
\begin{equation}\label{feh}
p^b_{-k,-k-\ell} \gamma_{k+\ell}=2^\ell\gamma_k, \quad
p^b_{-k,-k-\ell} \gamma_{k+\ell-1}=2^{\ell-1}\gamma_k.
\end{equation}
Applying (\ref{feh}) with $k=0$ we obtain using (\ref{ach1}) and
(\ref{ach2}) the equation
\begin{equation}\label{feh3}
n_0=2^{\ell-1}n^1_{-\ell}+2^\ell n^2_{-\ell}=
2^{\ell-1}\big(2n^1_{-\ell}+n^2_{-\ell}\big), \quad \ell \in
\mathbb{N}.
\end{equation}
Since (\ref{feh3}) holds for any $\ell \in \mathbb{N}$ but a nonzero
integer is not divisible by an arbitrary high power of $2$ we
conclude from (\ref{feh3}) that
\begin{equation}\label{feh4}
n_0=0.
\end{equation}
Applying (\ref{feh}) for $k \in \mathbb{N}$ and again using
(\ref{ach1}) and (\ref{ach2}) we get the equation
\begin{equation}\label{sto}
n^1_{-k}=0, \quad n^2_{-k}=2^{\ell-1}\big(2 n^1_{-k-\ell}+
n^2_{-k-\ell}\big),\quad k,\ell \in \mathbb{N}.
\end{equation}
The same reasoning which was used in the derivation of (\ref{feh4})
leads now to
\begin{equation}\label{sto2}
n^2_{-k}=0, \quad k \in \mathbb{N}.
\end{equation}
Hence the above three formulas give (\ref{ups}) and therefore
(\ref{feh2}). We conclude that
$$\overline{HM}=\lim_{\substack{\longrightarrow\\ b \to
\infty}}\lim_{\substack{\longleftarrow\\
a \to -\infty}}HM_a^b=0.$$
This finishes the proof of Step\,2.
\\ \\
\textbf{Step\,3: } \emph{$HM \neq 0$.}
\\ \\
Choose a sequence $\{a_j\}_{j \in \mathbb{N}}$ with
$a_j \in \mathbb{Z}$ for all $j \in \mathbb{N}$ which
satisfies the following two conditions.
\begin{itemize}
 \item $\lim_{k \to \infty} \sum_{j=1}^k2^{j-1}a_j=\infty$,
 \item $0< \frac{1}{2^k}\sum_{j=1}^k 2^{j-1}a_j<3/4$ for all $k$.
\end{itemize}
(Such a sequence can be easily constructed consisting only of zeroes
and ones).
We show that the element
$$\xi=\sum_{j=1}^\infty a_j \underline{c}_j$$
gives rise to a nonvanishing class in $HM$. Obviously
$\xi$ is in the kernel of the boundary operator. To show
that it is not in the image we argue by contradiction and
assume that there exists
$$\eta=\sum_{j=0}^\infty b_j \overline{c}_j$$
which coefficients $b_j \in \mathbb{Z}$ such that
$$\partial \eta=\xi.$$
It follows that
$$a_j=-2b_j+b_{j-1} ,\quad j \in \mathbb{N}.$$
By induction on this formula we obtain
$$b_0=\sum_{j=1}^k 2^{j-1}a_j+2^k b_k$$
for each $k \in \mathbb{N}$.
By our first assumption on the sequence $a_j$ we can find $n \in \mathbb{N}$
such that
$$\sum_{j=1}^k2^{j-1}a_j>b_0$$
for all $k\geq n$. The second assumption on the $a_j$ then implies
that
$$b_k = 2^{-k}b_0-2^{-k}\sum_{j=1}^k2^{j-1}a_j\in (-1,0)$$
for $k\geq n$ sufficiently large.
But $b_k$ is an integer. This contradiction shows that
$\xi$ does not lie in the image of $\partial$. This implies
Step\,3 and hence the theorem. \hfill $\square$

\bigskip
Kai Cieliebak, Ludwig-Maximilians-Universit\"at, 80333 M\"unchen, Germany\\
E-mail: kai@math.lmu.de
\\ \\
Urs Frauenfelder, Department of Mathematics and Research Institute
of
Mathematics, Seoul National University\\
E-mail: frauenf@snu.ac.kr

\end{document}